\documentclass[12pt]{amsart}

\usepackage{amssymb}
\usepackage{epsf}
\usepackage{amsmath,epsf}
\usepackage[latin1]{inputenc}
\usepackage{amsfonts}

\numberwithin{equation}{section}
\newtheorem{theo}{Theorem}[section]
\newtheorem{prop}[theo]{Proposition}

\newtheorem{lemm}[theo]{Lemma}

\newtheorem{rema}[theo]{Remark}

\newcommand{\pref}[1]{(\ref{#1})}

\def\cf{{\it cf. }}
\def\ie{{\it i.e. }}

\def\gex{\geq_{\rm lex}}

\def\gx{>_{\rm lex}}
\def\lx{<_{\rm lex}}

\def\N{{\mathbb N}}

\def\Q{{\mathbb Q}}

\def\I{{\mathcal I}}
\def\J{{\mathcal J}}
\def\A{{\mathcal A}}
\def\L{{\mathcal L}}
\def\D{{\mathcal D}}
\def\S{{\mathcal S}}

\def\G{{\mathcal G}}

\def\B{{\mathcal B}}

\def\P{{\mathcal P}}

\def\Qsym{{Qsym}}

\def\H{{\bf H}}
\def\SH{{\bf SH}}

\def\R{{\bf R}}

\def\shuffle{{\,\raise
1pt\hbox{$\scriptscriptstyle\cup{\mskip-4mu}\cup$}\,}}
\def\sshuffle{{\,\hbox{$\scriptscriptstyle\cup{\mskip-4mu}\cup$}\,}}

\title[Super-Covariant Polynomials]{Ideals of Quasi-Symmetric Functions and
Super-Covariant Polynomials for $\S_n$}

\author{J.-C.~Aval}
\address[J.-C. Aval]{Laboratoire A2X\\ Universit\'e Bordeaux 1\\ 351 cours
de 
 la Lib\'eration\\ 33405 Talence cedex\\ FRANCE}
\email{aval@math.u-bordeaux.fr}

\author{F.~Bergeron}
\address[F. Bergeron]{D\'epartement de Math\'ematiques\\ Universit\'e
 du Qu\'ebec \`a Montr\'eal\\ Montr\'eal, Qu\'ebec, H3C 3P8, CANADA}
\email{bergeron.francois@uqam.ca}

\author{N.~Bergeron}
 \address[N. Bergeron]{Department of Mathematics and Statistics\\ York
 University\\ 
     To\-ron\-to, Ontario M3J 1P3\\ CANADA}
\email{bergeron@mathstat.yorku.ca}
 \urladdr{http://www.math.yorku.ca/bergeron}

\thanks{F. Bergeron is supported in part by NSERC and
 FCAR}  
\thanks{N. Bergeron is supported in part by CRC, NSERC and PREA}

\date{\today} 

\begin{document} 

\begin{abstract} 
The aim of this work is to study the quotient ring $\R_n$ of the ring
$\Q[x_1,\dots,x_n]$ over the ideal $\J_n$ generated by non-constant
homogeneous quasi-symmetric functions. This article is a sequel to
\cite{a9}, in which is investigated the case of infinitely many variables.
We prove here that the dimension of $\R_n$ is given by $C_n$, the $n^{\rm
th}$ Catalan number. This is also the dimension of the space $\SH_n$ of
super-covariant polynomials, that is defined as the orthogonal complement of
$\J_n$ with respect to a given scalar product. We construct a basis for
$\R_n$ whose elements are naturally indexed by Dyck paths. This allows us to
understand the Hilbert series of $\SH_n$ in terms of number of Dyck paths
with a given number of factors.
\end{abstract} 

\maketitle

\section{Introduction}

We study, in this paper, a natural analog of the space $\H_n$ of covariant
polynomials of $\S_n$. Letting $\I_n$ denote the ideal generated by all
symmetric polynomials with no constant term
$$\I_n=\langle h_k,\ k>0\rangle,$$
where $h_k$ is the $k^{\rm th}$ homogeneous symmetric polynomials (\cf
\cite{mac}), the space $\H_n$ is defined as the orthogonal complement,
$\I_n^\perp$, in $\Q[x_1,\ldots,x_n]$, of the ideal $\I_n$, where the scalar
product considered is
\begin{equation}
\langle P,Q\rangle=P(\partial)Q(X)\big|_{X=0}.
\end{equation}
where $X$ stands for the variables $x_1,\ldots,x_n$,  $\partial$ stands for
$\partial x_1,\ldots,\partial x_1$, and in the same spirit, $X=0$ stands for
$x_1=\cdots=x_n=0$.

Equivalently (\cf \cite{orbit}, Proposition I.2.3), covariant polynomials
(also known as $\S_n$-harmonic polynomials) can be defined as polynomials
$P$ such that $Q(\partial)P=0$, for any symmetric polynomial $Q$ with no
constant term. Since, in particular, elements of $\H_n$ satisfy the Laplace
equation
$$(\partial x_1^2+\cdots+\partial x_n^2)\,P=\Delta\,P=0 ,$$
every covariant polynomial is also harmonic.

Classical results \cite{artin,steinberg} state that the space  $\H_n$
affords a graded
$\S_n$-module structure and is isomorphic (as a representation of $\S_n$) to
the left regular representation. Furthermore, as a graded $\S_n$-module,
$\H_n$ is isomorphic to the quotient
$$Q_n=\Q[X]/\I_n,$$
with $X=(x_1,\dots,x_n)$. The space $Q_n$ appears naturally in other
contexts; for
instance, as the cohomology ring of the
variety of complete flags \cite{borel}.
In particular, this implies that
\begin{equation}
\dim \H_n=n!\,. 
\end{equation}
Part of the interesting results surrounding the study of  $\H_n$ involve the
fact that it can also be described as the
linear span of all partial derivatives of the Vandermonde determinant. This
is just a special case of a more general result for finite
groups generated by reflections \cite{steinberg}.

By analogy, we consider here the space $\SH_n=\J_n^\perp$ of {\sl
super-covariant} polynomials, where
$\J_n$ is the ideal generated by {\sl quasi-symmetric} polynomials with no
constant term.
Since the ring of symmetric polynomials is a
subring of the ring of quasi-symmetric polynomials, we have $\I_n\subseteq
\J_n$  hence $\J_n^\perp\subseteq
\I_n^\perp$, thus
$$\SH_n\subseteq \H_n\,,$$
which somewhat justifies the terminology.
Quasi-symmetric polynomials where introduced by Gessel in 1984 \cite{ges}
and
have since appeared as a crucial tool in many interesting
algebraico-combinatorial contexts (\cf \cite{BMSW,NC,MR,stanley1,stanley2}).

As in the corresponding symmetric setup, we have a graded isomorphism
\begin{equation}
\SH_n\simeq \R_n=\Q[X]/\J_n
\end{equation}
and the approach used in the following work concentrates on this alternate
description. We will construct a basis of
$\R_n$ by giving an explicit set of monomial representatives. As we will
show, this set is naturally
indexed by {\sl Dyck paths} of length $n$, hence we obtain the following
main theorem.
\begin{theo}\label{main}
The dimension of $\SH_n$ is given by the well known Catalan numbers:
\begin{equation}
\dim\SH_n=\dim\R_n=C_n={1\over n+1}{2\,n\choose n}\,.
\end{equation}
In fact, taking into account the grading (with respect to degree), we have
the Hilbert series
\begin{equation}
\sum_{k=0}^{n-1}\dim \SH_n^{(k)}t^k=\sum_{k=0}^{n-1}\frac{n-k}{n+k}{n+k
\choose k} t^k\,.
\end{equation}
\end{theo}

The article is composed of five sections. In Section 2 we recall useful
definitions and basic properties. In Section 3 we construct a family $\G$ of
generators for the ideal $\J_n$ and state useful properties of this set. The
Section 4 is devoted to the proof of the main Theorem \ref{main}. We
construct an explicit basis for $\R_n$ which allows us in Section 5 to
obtain the Hilbert series of $\SH_n$.

\section{Basic definitions}

A {\sl composition} $\alpha =(\alpha_1,\alpha_2, \dots ,\alpha_k)$ of a
positive integer $d$ is an ordered list of positive integers ($>0$) whose
sum is $d$. We denote this by $\alpha \models d$  and also say that $\alpha$
is a composition of {\sl size} $d$. The size of $\alpha$ is denoted
$|\alpha|$. The integers
$\alpha_i$ are the {\sl parts} of $\alpha$, and the length $\ell(\alpha)$ is
set to be the
number of
parts of $\alpha$. 

There is a natural one-to-one correspondence between compositions of $d$ and
subsets of
$\{1,2,\ldots,d-1\}$. Let $S=\{ a_1,a_2,\dots ,a_k\}$ be such a subset, with
$a_1< \cdots <a_k$, then the
composition associated to $S$ is
$\alpha_d(S) =(a_1-a_0,a_2-a_1,\ldots ,a_{k+1}-a_k)$, where we set $a_0:=0$
and
$a_{k+1}:=d$. We denote $D(\alpha)$ the set associated to $\alpha$ by
this correspondence. For compositions $\alpha$ and $\beta$, we say that
$\beta$ is a {\sl refinement} of
$\alpha$,  if $D(\alpha)\subset D(\beta)$, and denote this by
$\beta\succeq \alpha$.

We will use vectorial notation for monomials. More precisely, for
$\nu=(\nu_1,\dots,\nu_n)\in\N^n$, we denote $X^\nu$ the monomial
$x_1^{\nu_1} x_2^{\nu_2} \cdots x_n^{\nu_n} $. We further denote
  $$[X^\nu]\,P(X)$$
the coefficient of the monomial $X^\nu$ in $P(X)$.

For a {\sl vector} $\nu\in\N^n$, let $c(\nu)$ the composition obtained by
erasing zeros (if any) in $\nu$. A polynomial $P\in\Q[X]$ is said to be {\sl
quasi-symmetric} if and only if, for any $\nu$ and $\mu$ in $\N^n$, we have
$$[X^\nu]P(X)=[X^\mu]P(X)$$
whenever $c(\nu)=c(\mu)$. The space of quasi-symmetric polynomials in $n$
variables is denoted
by $\Qsym_n$.
The space $\Qsym^{(d)}_n$ of homogeneous quasi-symmetric polynomials of
degree
$d$ admits as linear basis the set of {\sl monomial} quasi-symmetric
polynomials indexed by compositions of $d$. More precisely, for each
composition
$\alpha$ of $d$ with at most $n$ parts, we set
\begin{equation}
M_\alpha=\sum_{c(\nu)=\alpha} X^\nu
\end{equation}
For the $0$ composition, we set $M_0=1$. Another important linear basis is
that of the {\sl
fundamental} quasi-symmetric polynomials (\cf \cite{ges}):
\begin{equation}
F_\alpha=\sum_{\beta\succeq\alpha}M_\beta
\end{equation}
with $\alpha\models n$ and $\ell(\alpha)\leq n$. For example, with $n=4$,
\begin{eqnarray*} 
F_{21}(x_1,x_2,x_3,x_4)&=&M_{21}(x_1,x_2,x_3,x_4)+M_{111}(x_1,x_2,x_3,x_4)\\
&=&{x_1}^2\,x_2+{x_1}^2\,x_3+{x_1}^2\,x_4+{x_2}^2\,x_3+{x_2}^2\,x_4+{x_3}^2\
,x_4\\
&&\qquad +x_1\,x_2\,x_3+x_1\,x_2\,x_4+x_1\,x_3\,x_4+x_2\,x_3\,x_4.
\end{eqnarray*}

Part of the interest of fundamental quasi-symmetric functions comes from the
following properties. The first is trivial, but very useful and the second
comes from the theory of $P$-partitions \cite{stanley1,stanley2}.
\begin{prop}
For $\alpha =(\alpha_1,\alpha_2,\dots,\alpha_k)\models d$,
\begin{equation}\label{Frel}
F_\alpha(X) =
\left\{\begin{array}{ll}
x_1 F_{(\alpha_1-1,\alpha_2,\ldots,\alpha_k)}(X)  + F_\alpha(x_2,\dots,x_n)
&\mbox{ if $\alpha_1>1$,}\\
&\mbox{ }\\
x_1 
F_{(\alpha_2,\alpha_3,\ldots,\alpha_{k})}(x_2,\dots,x_n)+F_\alpha(x_2,\dots,
x_n)& \mbox{ if $\alpha_1=1$.}
\end{array}\right. 
\end{equation} 
\end{prop}

Let $u=u_1\cdots u_l\in\S_\ell$ and $v=v_1\cdots
v_m\in\S_{[\ell+1,\ell+m]}$. Let $u\shuffle v$ denote the set of {\sl
shuffles} of the words $u$ and $v$, \ie $u\shuffle v$ is the set of all
permutations $w$ of $\ell+m$ such that $u$ and $v$ are subwords of $w$. In
particular $u\shuffle v$ contains ${\ell+m \choose m}$ permutations. Let
$\D(u)=\{i,\ u_{i}>u_{i+1}\}$ denote the {\sl descent set} of $u$. If
$\beta$ and $\gamma$ are the two compositions such that $D(\beta)=\D(u)$ and
$D(\gamma)=\D(v)$, then

\begin{prop}[\cite{stanley2}, Exercise 7.93]\label{propprod}
\begin{equation}
F_\beta\,F_\gamma=\sum_{w\in u\sshuffle v} F_{\alpha_{\ell+m}(D(w))}.
\end{equation}
\end{prop}

There is an evident bijection between elements $\nu$ of $\N^n$ and
the corresponding monomial $X^\nu$.
Elements of $\N^n$ are naturally called vectors. Just as for compositions,
the size $\nu_1+\cdots+\nu_n$ of $\nu$ is denoted $|\nu|$. It will also be
convenient to denote $\ell(\nu)$ the position of its last non-zero
component. As usual,  $\nu+\mu$ is the componentwise addition of vectors.

To make for easier reading, we generally use $\alpha,\ \beta,\ \gamma$ to
denote compositions, and $\mu,\ \nu$ to denote vectors. In general, $n$ the
length of vectors (or number of variables) is fixed, and if $w$ is a word of
integers  (that is an element of $\N^k$ for $0\le k\le n$) we denote by
$w0^*=w0^{n-k}$ the vector whose first $k$ parts are the {\sl letters} of
$w$, to which are added $n-k$ $0$'s at the end. If $u=u_1\cdots u_k$ and
$v=v_1\cdots v_m$ are words of integers,  the word
  $$u\,v:=u_1\cdots u_k v_1\cdots v_m$$
is the {\sl concatenation} of $u$ and $v$. We use the same symbol $\alpha$
for both the composition $(\alpha_1,\dots,\alpha_\ell)$ and the word
$\alpha_1,\dots,\alpha_\ell$, likewise for vectors.

We next associate to any vector $\nu$ a path $\pi(\nu)$ in the $\N\times\N$
plane with steps going north or east as follows. If
$\nu=(\nu_1,\dots,\nu_n)$, the path $\pi(\nu)$ is
$$(0,0)\rightarrow(\nu_1,0)\rightarrow(\nu_1,1)\rightarrow(\nu_1+\nu_2,1)
\rightarrow(\nu_1+\nu_2,2)\rightarrow\cdots\ \ \ \ \ \ \ \ \ \ \ \ \ \ \ \ \
\
\ \ \ $$
\vskip -0.6cm
$$\ \ \ \ \ \ \ \ \ \ \ \ \ \ \ \ \ \ \ \ \ \ \ \ \ \ \ \ \ \ \ \ \ \ \ \ \
\rightarrow(\nu_1+\cdots+\nu_n,n-1)\rightarrow(\nu_1+\cdots+\nu_n,n).$$
For example the path associated to $\nu=(2,1,0,3,0,1)$ is

\vskip 0.2cm
\centerline{\epsffile{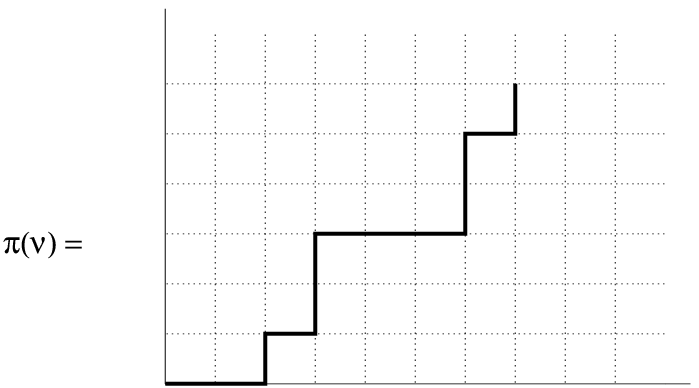}}

\vskip 0.2cm
\noindent Observe that the height of the path is always $n$, whereas its
{\sl width} is $|\nu|$.

We distinguish two kinds of paths, thus two kinds of vectors, with respect
to their ``behavior'' regarding the diagonal $y=x$.
If the path remains above the diagonal, we call it a {\sl Dyck path}, and
say that the corresponding vector is {\sl Dyck}. If not, we say that the
path (or equivalently the associated vector) is {\sl transdiagonal}. For
example $\eta=(0,0,1,2,0,1)$ is Dyck and $\varepsilon=(0,3,1,1,0,2)$ is
transdiagonal.

\vskip 0.2cm
\centerline{\epsffile{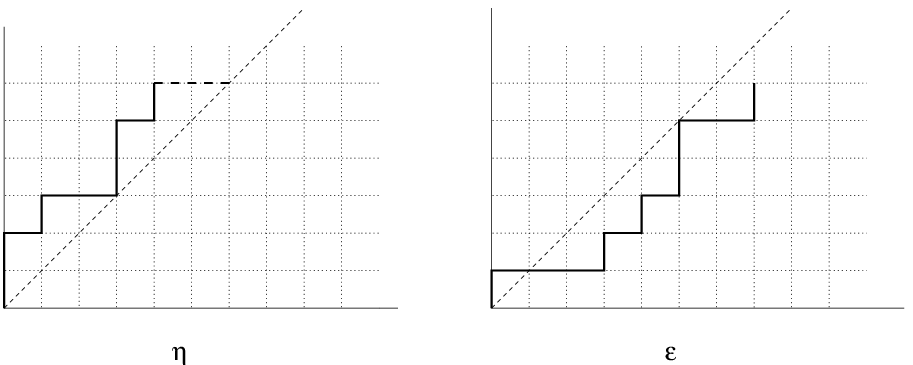}}

\vskip 0.2cm
Observe that $\nu=\nu_1\cdots\nu_n$ is transdiagonal if and only if there
exists $1\le \ell\le n$ such that
\begin{equation}\label{cross}
\ell<\nu_1+\ldots+\nu_\ell.
\end{equation}

Recall that the classical lexicographic order, on monomials of same degree,
is
\begin{equation}
X^\nu\gex X^\mu\qquad{\rm iff}\qquad \nu\gex \mu,
\end{equation}
where we say that $\nu$ is lexicographically larger than $\mu$,
$\nu\gx\mu$, if the first non-zero part of the
vector $\nu-\mu$ is positive. Thus
  $$x_1^3 \gx x_1^2x_2\gx x_1x_2^2\gx x_2^3$$
since
  $$(3,0)\gx(2,1)\gx(1,2)\gx(0,3)\,.$$
We extend this order to all monomials (of possibly different degree) by
setting
  $$X^\nu\lx X^\mu\qquad{\rm whenever}\qquad |\nu|<|\mu|\,.$$
This is known as the {\sl graded lex} order, and it clearly makes sense for
vectors.  

\section{The $\G$ basis}

Following \cite{a9}, we exploit relations \pref{Frel} to construct a family
$$\G=\{G_\varepsilon\} \subset\J_n$$
indexed by vectors that are transdiagonal. For $\alpha$ any composition of
$k\leq n$, the  polynomial $G_\varepsilon$, with $\varepsilon:=\alpha0^*$,
is defined to be 
\begin{equation}
G_\varepsilon:=F_\alpha.
\end{equation}
When $\alpha\not=0$, the vector $\varepsilon=\alpha0^*$ is clearly
transdiagonal. For a general vector $\varepsilon$ (not of the form
$\alpha0^*$), the  polynomial $G_\varepsilon$ is defined recursively in the
following way. Let $\varepsilon=w0a\beta0^*$ be the unique factorization of
$\varepsilon$ such that $w$ is a word of $k-1$ non-negative integers, $a>0$
is a positive integer, and $\beta$ is a composition (parts $>0$). Then we
set
\begin{equation}\label{Geps}
G_\varepsilon=G_{wa\beta0^*}-x_k\,G_{w(a-1)\beta0^*}.
\end{equation}
Both terms on the right of \pref{Geps} are well defined, and moreover we
have
\begin{itemize}
\item $\ell(wa\beta0^*)=\ell(w(a-1)\beta0^*)=\ell(\varepsilon)-1$;
\item  $wa\beta0^*$ and $w(a-1)\beta0^*$ are transdiagonal as soon as
$\varepsilon$ is transdiagonal.
\end{itemize}
In fact, let $\ell$ be the first ordinate where $\pi(\varepsilon)$ crosses
the diagonal, this is to say that it is the smallest integer such that
$\ell<\varepsilon_1+\ldots+\varepsilon_\ell$. Then the second assertion
follows from
$$\varphi_1+\ldots+\varphi_\ell>\psi_1+\ldots+\psi_\ell
      =\varepsilon_1+\ldots+\varepsilon_\ell-1>\ell-1,$$
where $\varphi=wa\beta0^*$ and $\psi=w(a-1)\beta0^*$.

For example,
\begin{eqnarray*}
G_{1020}&=&G_{1200}-x_2\,G_{1100}\\
&=&F_{12}(x_1,x_2,x_3,x_4)-x_2\,F_{11}(x_1,x_2,x_3,x_4)\cr
&=&x_1\,{x_2}^2+x_1\,{x_3}^2+x_1\,{x_4}^2+x_2\,{x_3}^2+x_2\,{x_4}^2+x_3\,{x_
4}^2\\
&&\qquad\qquad +x_1\,x_2\,x_3+x_1\,x_2\,x_4+x_1\,x_3\,x_4+x_2\,x_3\,x_4\\
&&\qquad-x_2\,(x_1\,x_2+x_1\,x_3+x_1\,x_4+x_2\,x_3+x_2\,x_4+x_3\,x_4)\\
&=&x_1\,{x_3}^2+x_1\,x_3\,x_4+x_1\,{x_4}^2-{x_2}^2\,x_3-{x_2}^2\,x_4+x_2\,{x
_3}^2+x_2\,{x_4}^2+x_3\,{x_4}^2.\\
\end{eqnarray*}
Observe on this example that the leading monomial (in graded lex order) of
$G_{1020}$ is $X^{1020}=x_1^1x_2^0x_3^2x_4^0$.
This holds in general for the $\G$ family as stated in the following
proposition, for which all technical details can be found in \cite{a9}.
\begin{prop}[\cite{a9}, Corollary 3.4]
The leading monomial $LM(G_\varepsilon)$ of $G_\varepsilon$ is
$X^\varepsilon$.
\end{prop}

\section{Proof of the main theorem}

We now prove our main Theorem \ref{main}, and more precisely obtain an
explicit basis for the space $\R_n$ naturally indexed by Dyck paths, thus of
cardinality equal to $C_n$.
\begin{theo}\label{main2}
The set of monomials
\begin{equation}
\B_n=\{X^\eta\ |\ \pi(\eta)\ {\it is\ a\ Dyck\ path}\}
\end{equation}
is a basis of the space $\R_n$.
\end{theo}
The proof will be achieved in several steps. We start with the following
lemma.
\begin{lemm}\label{bsup}
Any polynomial $\P\in\Q[X]$ is in the linear span of $\B_n$ modulo $\J_n$,
which is to say that
\begin{equation}\label{redeq}
P(X)\equiv\sum_{X^\eta\in\B_n} c_\eta X^\eta\ \ \ \ ({\rm mod}\ \J_n).
\end{equation}
\end{lemm}
\proof
It clearly suffices to show that \pref{redeq} holds for any monomial
$X^\nu$,
with $\nu$ transdiagonal. We assume that there exists $X^\nu$ not reducible
in the form \pref{redeq} and we choose $X^\varepsilon$ to be the smallest
amongst them with respect to the lexicographic order. Let us write
\begin{eqnarray*}
X^\varepsilon&=&LM(G_\varepsilon)\\
&=&(X^\varepsilon-G_\varepsilon)+G_\varepsilon\\
&\equiv&X^\varepsilon-G_\varepsilon\ \ \ \ ({\rm mod}\ \J_n).
\end{eqnarray*}
All monomial in $(X^\varepsilon-G_\varepsilon)$ are smaller than
$X^\varepsilon$, thus they are reducible. This contradicts our assuption on
$X^\varepsilon$ and completes our proof.
\endproof

Thus $\B_n$ spans the space $\R_n$. We now prove its linear independence.
This is equivalent to showing that the set $\G$ is a Gr\"obner basis of the
ideal $\J_n$. A crucial lemma is the following one, which is the
quasi-symmetric analogue of a classical result is the case of symmetric
polynomials (\cite{orbit}, Theorem II.2.2).
\begin{lemm}\label{lele}
If we denote by $\L[S]$ the linear span of a set $S$, then
\begin{equation}\label{lmp}
\Q[X]=\L[X^\eta F_\alpha,\ X^\eta\in\B_n,\ \alpha\models r\ge 0].
\end{equation}
\end{lemm}
\proof
We have already obtained
$$X^\varepsilon=\sum_{X^\eta\in\B_n} c_\eta X^\eta\ \ \ \ ({\rm mod}\
\J_n),$$
which is equivalent to
\begin{equation}\label{la}
X^\varepsilon=\sum_{X^\eta\in\B_n} c_\eta X^\eta + \sum_{\alpha\models r\ge
1} Q_\alpha F_\alpha.
\end{equation}
We then apply the reduction \pref{la} to each monomial of the $Q_\alpha$'s
and use Proposition \ref{propprod} to reduce products of fundamental
quasi-symmetric functions. We obtain \pref{lmp} in a finite number of
operations since degrees strictly decreases at each operation, because
$\alpha\models r\ge 1$ implies $\deg F_\alpha\ge 1$.
\endproof

The next lemma is the final step in our proof of the Theorem \ref{main2}.
\begin{lemm}\label{crux}
The set $\G$ is a linear basis of the ideal $\J_n$, \ie
\begin{equation}
\J_n=\L[G_\varepsilon\ | \ \varepsilon\ {\it transdiagonal}].
\end{equation}
\end{lemm}
\proof
Let us denote by $\A_n$ the set
\begin{equation}
\A_n=\{X^\xi\ |\  x_1^{\xi_n}\,x_2^{\xi_{n-1}}\cdots x_n^{\xi_1}\in \B_n\}
\end{equation}
Now the algebra endomorphism of $\Q[X]$ that {\sl reverses} the variables,
that is
   $$x_i\mapsto x_{n-i+1},$$
clearly fixes the subalgebra $\Qsym$. In fact it maps $F_\alpha$ to
$F_{\alpha'}$, where $\alpha'$ is the reverse composition.

It follows from Lemma \ref{lele} that:
\begin{equation}\label{eqjolie}
\Q[X]=\L[X^\xi F_\alpha\ | \ X^\xi\in\A_n,\ \alpha\models r\ge 0].
\end{equation}
Now to prove Lemma \ref{crux}, we reduce the problem as follows. We first
use \pref{eqjolie} and Proposition \ref{propprod} to write
$$\begin{array}{ll}
\J_n&\!\!\!=\langle F_\alpha,\ \alpha\models s\ge
0\rangle_{\Q[X]}=\L[X^\xi\, F_\alpha\, F_\beta\ | \ X^\xi\in\A_n,\
\alpha\models s\ge 0,\ \beta\models t\ge 1]\\
&\mbox{ }\\
&\!\!\!=\L[X^\xi\, F_\gamma\ |\ X^\xi\in\A_n,\ \gamma\models r\ge 1].
\end{array}$$
It is now sufficient to prove that for all $X^\xi\in\A_n$ and all $
\gamma\models r\ge 1$
\begin{equation}\label{eqfin}
 X^\xi\,F_\alpha\in\L[G_\varepsilon\ | \ \varepsilon\ {\it transdiagonal}].
\end{equation}
But Lemma \ref{bsup} implies that any monomial of degree greater than $n$ is
in
$\J_n$. Hence to prove \pref{eqfin}, we need only show it for $\xi$ and
$\gamma$ such that
$|\xi|+|\gamma|\le n$. To do that, we reduce the product
\begin{equation}\label{rea}
 x_{n}^{\eta_{n}}(x_{n-1}^{\eta_{n-1}}(\cdots(x_2^{\eta_2}(x_1^{\eta_1}
F_\alpha))))
\end{equation}
recursively, using
\begin{equation}\label{re1}
x_k\,G_{wb\beta 0^*}=G_{w(b+1)\beta 0^*}-G_{w0 (b+1)\beta 0^*}
\end{equation}
or
\begin{equation}\label{re2}
x_k\,G_{w0^*00^*}=G_{w0^*10^*}-G_{w0^*010^*}.
\end{equation}
Relations \pref{re1} and \pref{re2} are immediate consequences of the
definition of
the $\G$ basis (relation \pref{Geps}).

We have to show that the vectors $\varepsilon$ generated in this process are
all transdiagonal and that the length $\ell(\varepsilon)$ always remains at
most equal to $n$. 
Let us first check that the transdiagonal part. This is obvious in the
case of relation \pref{re2}. In the other case (relation \pref{re1}), it is
sufficient to observe that, for $\varphi=wb\beta 0^*$, if $m$ is such that
   $$\varphi_1+\ldots+\varphi_m>m$$
with $m>\ell(w)$ (if not, it is evident), then
  $$\varphi'_1+\ldots+\varphi'_m>m+1>m \qquad {\rm and}\qquad
\varphi_1+\ldots+\varphi_{m+1}>m+1.$$
where $\varphi'=w(b+1)\beta 0^*$, and $\varphi''=w0(b+1)\beta 0^*$.
We shall now prove that the length of the $\varepsilon$'s always remains at
most equal to $n$.
For this we need to keep track of the term
$\varepsilon_{\ell(\varepsilon)}$. Two cases have to be considered.
\begin{itemize}
\item First case: $\varepsilon_{\ell(\varepsilon)}$ comes from
$\alpha_{\ell(\alpha)}$ that has been shifted on the right by relation
\pref{re2}. It could have made at most $|\xi|$ steps on the right, whence
$$\ell(\varepsilon)\le \ell(\alpha)+|\xi|\le |\alpha|+|\xi|\le n.$$
\item Second case: $\varepsilon_{\ell(\varepsilon)}$ is a $1$ generated by
relation \pref{re2} that has been shifted on the right. If it is generated
by
multiplication by $x_k$, we consider the vector
$$\eta=\xi_n\xi_{n-1}\cdots\xi_{k}0^*.$$
Since $X^\xi\in\A_n$ implies $\pi(\eta)$ is a Dyck path, we have
$$|\eta|<\ell(\eta)=n-k+1$$
hence the $1$ generated can be shifted at most in position
$$k+|\eta|\le k+n-k=n.$$
\end{itemize}
\endproof

 The recursive process used to reduced a product of form  \pref{rea}  is
illustrated in the following example, where $n=5$.
$$\begin{array}{ll}
x_1\,x_3\,F_{21}&\!\!\!=x_3(x_1\,F_{21})\\
&\mbox{ }\\
&\!\!\!=x_3(G_{31000}-G_{03100})\\
&\mbox{ }\\
&\!\!\!=x_3\,G_{31000}-x_3\,G_{03100}\\
&\mbox{ }\\
&\!\!\!=G_{31100}-G_{31010}-G_{03200}+G_{03020}.
\end{array}$$

\proof[End of proof of Theorem \ref{main2}]
By Lemma \ref{bsup}, the set $\B_n$ spans the quotient $\R_n$, and we are
now in a position to prove its linear independence. Assume we have a linear
dependence relation modulo $\J_n$, \ie there exists $P$
$$P=\sum_{X^\xi\in\B_n} a_\xi X^\xi\,\in\,\I_n.$$
By Lemma \ref{crux}, $\J_n$ is linearly spanned by the $G_\varepsilon$'s,
thus
$$P=\sum_{\varepsilon\ {\rm transdiagonal}} b_\varepsilon G_\varepsilon.$$
This implies $LM(P)=X^\varepsilon$, with $\varepsilon$ transdiagonal, which
is absurd.
\endproof

\section{Hilbert series}

Since Theorem \ref{main2} gives us an explicit basis for the quotient
$\R_n$, which is isomorphic to $\SH_n$ as a graded vector space, we are able
to refine Theorem \ref{main} by giving the Hilbert series of the space of
super-covariant polynomials. For $k\in\N$, let $\SH_n^{(k)}$ and
$\R_n^{(k)}$ denote the projections
\begin{equation}
\SH_n^{(k)}=\SH_n\,\cap\,\Q^{(k)}[X]\simeq
\R_n\,\cap\,\Q^{(k)}[X]=\R_n^{(k)}
\end{equation}
where $\Q^{(k)}[X]$ is the vector space of homogeneous polynomials of degree
$k$ together with zero.
Here, we represent Dyck path horizontally, with $n$ rising steps $(1,1)$
and $n$ falling steps $(1,-1)$. Let us denote by $D_n^{(k)}$ the number of
Dyck paths of length $2n$ ending by exactly $k$ falling steps and by
$C_n^{(k)}$ the number of Dyck paths of length $2n$ which have exactly $k$
factors, \ie $k+1$ points on the axis. The next figure gives an example of a
Dyck path of length $28$, ending with $4$ falling steps and made of $3$
factors.

\vskip 0.2 cm
\centerline{\epsffile{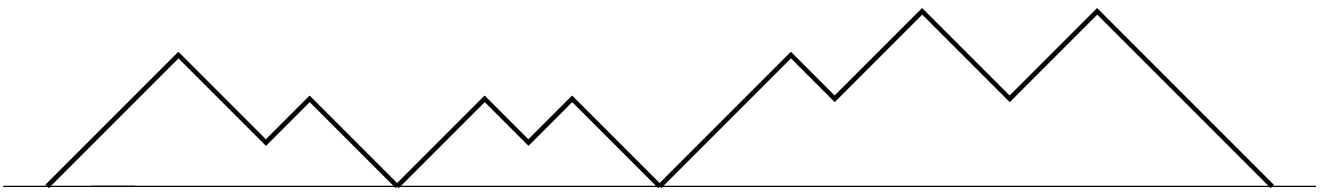}}

\vskip 0.2 cm
It is well known that
\begin{equation}
D_n^{(k)}=C_n^{(k)}=\frac{k\,(2n-k-1)!}{n!\,(n-k)!}\raise 2pt\hbox{,}
\end{equation}
where the first equality is classical (\cf \cite{vaille} for example for a
bijective proof), and the
second corresponds to \cite{Kreweras}, formula (7).

Let us denote by $F_n(t)$ the Hilbert series of $\SH_n$, \ie
\begin{equation}
F_n(t)=\sum_{k\ge 0} \dim \SH_n^{(k)}\, t^k.
\end{equation}

\begin{theo}
For $0\le k\le n-1$, the dimension of $\SH_n^{(k)}$ is given by
\begin{equation}
\dim \SH_n^{(k)}=\dim \R_n^{(k)}=D_n^{(n-k)}=C_n^{(n-k)}=\frac{n-k}{n+k}{n+k
\choose k}\,.
\end{equation}
For $k\ge n$ the dimension of $\SH_n^{(k)}$ is $0$.
\end{theo}
\proof
By Theorem \ref{main2}, we know that the set
$$\B_n=\{X^\eta\ |\ \pi(\eta)\ {\it is\ a\ Dyck\ path}\}$$
is a basis for $\R_n$. It is then sufficient to observe that the path
$\pi(\eta)$ associated to $\eta$ ends by exactly $n-|\eta|$ falling steps.
\endproof

For example, we have:

\medskip
\hskip 1cm
\begin{tabular}{|l|l|}
\hline
$n$ & $F_n(t)$ \\
\hline
1 & $1$ \\ 
\hline
2 & $1+t$ \\ 
\hline
3 & $1+2t+2t^2$ \\ 
\hline
4 & $1+3t+5t^2+5t^3$ \\
\hline
5 & $1+4t+9t^2+14t^3+14t^4$ \\
\hline
6 & $1+5t+14t^2+28t^3+42t^4+42t^5$ \\
\hline
7 & $1+6t+20t^2+48t^3+90t^4+132t^5+132t^6$ \\
\hline
\end{tabular}

\vskip 0.2 cm
This gives
 \begin{equation}
F_n(t)=\sum_{k=0}^{n-1}\frac{n-k}{n+k}{n+k \choose k} t^k\,.
\end{equation}
from which one easily deduces that the generating series for the $F_n(t)$'s
is
 \begin{equation}
\sum_n F_n(t)\,x^n={1-\sqrt{1-4\,t\,x}-2\,t\over 2\,(t+x-1)}.
\end{equation}

\begin{rema}\rm
The study of various filtrations of the space $\Q[X]$, with respect to
family of ideals of quasi-symmetric polynomials, will be the object of a
forthcoming paper \cite{suite}.
\end{rema}

\section*{Acknowledgments}

\noindent
We thank Adriano Garsia for stimulating discussions about this work.


\end{document}